\documentclass{amsart}
\usepackage{amscd} 
\usepackage{amsfonts} 
\usepackage{amssymb} 
\usepackage{latexsym} 
\input{diagrams}
\newcommand{\ncm}{\newcommand}

\newtheorem{theorem}{Theorem}[section]
\newtheorem{prop}[theorem]{Proposition}

\newtheorem{cor}[theorem]{Corollary}

\newtheorem{lem&def}[theorem]{Lemma \& Definition}
\newtheorem{definition}[theorem]{Definition}
\newtheorem{example}[theorem]{Example}

\ncm{\End}{\mbox{\rm End}\,}
\def\Aut{\mbox{\rm Aut}}
\def\Hom{\mbox{\rm Hom}\,}

\def\Im{\mbox{\rm Im}\,}

\def\id{\mbox{\rm id}}

\def\into{\hookrightarrow}
\def\to{\rightarrow}

\def\o{\otimes}    
\def\b{\bowtie}    

\ncm{\rarr}[1]{\stackrel{#1}{\longrightarrow}}
\ncm{\larr}[1]{\stackrel{#1}{\longleftarrow}}
\def\cop{\Delta}

\def\eps{\varepsilon}

\def\du1{\hat 1}

\def\-1{_{(-1)}}
\def\0{_{(0)}}
\def\1{_{(1)}}
\def\2{_{(2)}}
\def\3{_{(3)}}
\def\4{_{(4)}}

\def\|{\, | \,}

\def\du1{\hat 1}
\def\lact{\triangleright}

\begin{document}

\title[Pseudo-Galois Extensions
and Hopf Algebroids]{Pseudo-Galois Extensions
and  Hopf Algebroids}
\author{Lars Kadison}
\address{Department of Mathematics \\University of Pennsylvania \\ David Rittenhouse
Laboratory \\ 209 South 33rd Street \\ Philadelphia, PA 19104-6395} 
\email{lkadison@c2i.net} 
\date{}
\thanks{The author thanks Caenepeel, Van Oystaeyen
and Torrecillas for a pleasant stay in 
Belgium in September 2005.}

\begin{abstract} 
 A pseudo-Galois extension  is shown to be a depth two extension.
Studying its left bialgebroid, we construct an enveloping Hopf algebroid
for the semi-direct product of groups, or more generally involutive Hopf algebras, and their module algebras. It  is a type of cofibered
sum of two inclusions of the Hopf algebra into the semi-direct
product and its derived right crossed product.
Van Oystaeyen and Panaite observe that this Hopf algebroid is non-trivially isomorphic
to a Connes-Moscovici Hopf algebroid, which raises interesting
comparative questions.     
\end{abstract} 
\maketitle

\section{Introduction}

The analytic notion of finite depth for subfactors was widened to
the algebraic setting of Frobenius extensions
in \cite{KN}, the main theorem of which states that a certain depth two  Frobenius
extension  $A \| B$ with trivial centralizer is a Hopf-Galois extension.
The theorem and its proof is   essentially  a reconstruction theorem, which uses an Ocneanu-Szymanski  pairing of two centralizers on the tower of algebras
above $A \| B$, isomorphic to the two main players in this paper 
$\End {}_BA_B$ and $\End {}_A A \! \o_B A_A$; then shows that the resulting algebra-coalgebra is
a Hopf algebra.  In the paper \cite{KS} the notion of depth two was widened to
arbitrary algebra extensions whereby it was shown that the  bimodule endomorphism ring
$\End {}_BA_B$ of a  depth two extension $A \| B$
has a bialgebroid structure as in Lu \cite{Lu}. Interesting classes of examples were noted such as 
finite dimensional algebras, weak Hopf-Galois extensions, H-separable extensions
and various normal subobjects in quantum algebra.  Later in \cite{fer}, the underlying fact
emerged  that any (including infinite index) algebra extension $A \| B$ is right Galois w.r.t.\ a bialgebroid coaction if and only if it is right depth two and the natural module $A_B$ is balanced.  This is in several respects analogous
to the characterization of a Galois field extension as normal and separable (the splitting
field of separable polynomials).    

A  bialgebroid
 is in simplest terms a bialgebra over a noncommutative base,
which are of interest from the point of view of tensor categories \cite{ENO}
and mathematical physics \cite{BB}.  
A depth two extension $A \| B$  has a Galois action machinery consisting 
of a bialgebroid with base algebra the centralizer $C_A(B)$ of the extension. For example,
if $A$ is a commutative ring, then this Galois bialgebroid specializes to the Galois biring
in Winter \cite{W}.  If $A \| B$ is Frobenius and $C_A(B)$ is semisimple in all base field extensions,
as  with reducible subfactors, the extension is a weak Hopf-Galois extension,
where the Galois bialgebroid is a weak Hopf algebra \cite{KN, CDG}.  
 Bialgebroids equipped with antipodes are Hopf algebroids, a recent object of study with competing definitions of what constitutes
an antipode \cite{BS, Lu, KR}.
   In Section~3, we investigate the Galois bialgebroid of a new type of depth two extension called  a pseudo-Galois extension,  
which is a notion  generalizing the notions of H-separable extension (e.g.\ an Azumaya algebra) and group-Galois extension (e.g.\ a Galois algebra)
 \cite{MM}.
Its Galois bialgebroid is shown in Theorem~\ref{th-epic} to be closely related to a certain Hopf algebroid (in the sense of B\"ohm and Szlach\'anyi \cite{BS}) we obtain from a cofibered sum of a semi-direct product
of an algebra with a group of automorphisms and its opposite right crossed product.
This Hopf algebroid extends Lu's basic Hopf algebroid on
the enveloping  algebra over an algebra \cite{Lu}, 
to the group action setting, and 
in Theorem~\ref{th-bowtie} to the involutive Hopf algebra action setting. In the last section of this paper we discuss
Van Oystaeyen and Panaite's isomorphism of
this  ``enveloping'' Hopf algebroid given in the preprint \texttt{QA/0508411} to this  paper
 with a
Hopf algebroid 
 in Connes-Moscovici \cite{CM}, which is known as a para-Hopf algebroid in \cite{KR}. The
isomorphism  may be derived from a universal condition in Proposition~\ref{prop-universal}
or from the unit representation of a bialgebroid on its base algebra. 
The isomorphism of bialgebroid invariants naturally raises the  possibility of
a relation between certain pseudo-Galois extensions, or a generalization we propose,  with Rankin-Cohen brackets \cite{CM}.

\section{Depth two and pseudo-Galois extensions}

All algebras in this paper are unital associative algebras over a commutative ground ring $K$.
An algebra extension $A \| B$ is a unit-preserving algebra homomorphism $B \to A$,
a proper extension if this mapping is monic.  The induced
bimodule ${}_BA_B$ plays the main role below.  Unadorned tensors,  hom-groups and endomorphism-groups (in a category of modules) between algebras are
over the ground ring unless otherwise stated.
For example, $\End A$ denotes the linear
endomorphisms of an algebra $A$, but not the algebra
endomorphisms of $A$; and $\Hom (A_B, B_B)$
denotes the $B$-$A$-bimodule of right $B$-bimodule homomorphisms from $A$ into $B$. The default setting is the natural module
structure unless otherwise specified.     

An algebra extension $A \| B$ is \textit{left depth two (D2)} if
its tensor-square $A \o_B A$ as a natural $B$-$A$-bimodule
is isomorphic to a direct summand of a finite direct sum
of the natural $B$-$A$-bimodule $A$: for some positive integer
$N$, we have 
\begin{equation}
\label{eq: D2}
A \o_B A \oplus * \cong A^N
\end{equation}
An extension $A \| B$ is \textit{right D2} if eq.~(\ref{eq: D2}) holds instead as natural
$A$-$B$-bimodules.  

 Since condition~(\ref{eq: D2})
implies maps in  two hom-groups satisfying $\sum_{i=1}^N g_i \circ f_i = \id_{A \o_B A}$, where  $g_i \in \Hom ({}_BA_A, {}_BA\o_B A_A) \cong (A \o_B A)^B$ (via $g \mapsto g(1)$) and $$f_i \in \Hom ({}_BA\o_B A_A, {}_BA_A) \cong
\End {}_BA_B := \mathcal{S} $$
via $f \mapsto (a \mapsto f(a \o_B 1))$, we obtain an equivalent condition for extension $A \| B$ to be
left D2: there is a positive integer $N$, $\beta_1,\ldots,\beta_N \in \mathcal{S}$
and $t_1,\ldots,t_N \in (A \o_B A)^B$ (i.e., satisfying for each $i = 1,\ldots,n$, $bt_i = t_i b$ for every $b \in B$) such that 
\begin{equation}
\label{eq: lqbd2}
\sum_{i=1}^N t_i \beta_i(x)y = x \o_B y
\end{equation}
for all $x,y \in A$.

Like dual bases for projective modules,
this equation is  useful. For example, to show
$\mathcal{S}$ finite projective as a left $C_A(B)$-module (module action
given by $r \cdot \alpha = \lambda_r \circ \alpha$), apply 
$\alpha \in \mathcal{S}$ to the first tensorands of the equation, set $y = 1$ and apply the
multiplication mapping $\mu: A \o_B A \to A$ to obtain
\begin{equation}
\label{eq: tys}
 \alpha(x)  = \sum_i \alpha(t^1_i) t_i^2 \beta_i(x), 
\end{equation}
where we suppress a possible summation in $t_i \in A \o_B A$ using a Sweedler notation, 
$t_i = t_i^1 \o_B t^2_i$.
But for each $i = 1, \ldots,N$, we note that $$T_i(\alpha) := \alpha(t^1_i) t^2_i \in C_A(B) := R $$ 
defines a homomorphism $T_i \in \Hom ({}_R\mathcal{S}, {}_RR)$, so that
eq.~(\ref{eq: tys}) shows that $\{ T_i \}$, $\{ \beta_i \}$ are finite dual bases for ${}_R\mathcal{S}$.  

As another example, Eq.~(\ref{eq: lqbd2}) is used in the explicit formula for coproduct
in Eq.~(\ref{eq: so}), which in the case $A$ is a commutative ring
and $B$ a subfield would give an explicit formula for ``preservations'' in \cite{W}
in terms of dual bases and implying a simpler  proof for the Galois biring correspondence theorem \cite[Theorem 6.1]{W}. (We note that the generalized Jacobson-Bourbaki correspondence theorem \cite[Theorem 2.1]{W} may be extended with almost the same argument to any noncommutative algebra $A$ possessing algebra homomorphism into a division algebra.)

Similarly, an algebra extension $A \| B$ is right D2 if there is a positive
integer $N$, elements $\gamma_j \in \End {}_BA_B$ and
$u_j \in (A \o_B A)^B$ such that
\begin{equation}
\label{eq: right D2}
x \o_B y = \sum_{j=1}^N x \gamma_j(y)u_j
\end{equation}
for all $x,y \in A$.  We call the elements $\gamma_j \in \mathcal{S}$ and $u_j \in (A \o_B A)^B$
right D2 quasibases for the extension $A \| B$. Fix this notation and the corresponding notation
$\beta_i \in \mathcal{S}$ and $t_i \in (A \o_B A)^B$ for left D2 quasibases throughout this paper. 
An algebra extension is of course D2 if it is both left and right D2. 

Recall from \cite{KS} that a depth two extension $A \| B$ has left bialgebroid structure on
the algebra $\End {}_BA_B$ with noncommutative base algebra $C_A(B)$.  Again let $R$ denote $C_A(B)$ and
$\mathcal{S}$ denote the bimodule endomorphism algebra $\End {}_BA_B$  (under composition).  We sketch
the left $R$-bialgebroid structure on $\mathcal{S}$ since we will  need it in the case 
of pseudo-Galois extensions.  

Its $R$-$R$-bimodule structure ${}_R\mathcal{S}_R$ is generated by the algebra homomorphism $\lambda: R \to \mathcal{S}$ given
by $r \mapsto \lambda_r$, left multiplication by $r$, and the algebra anti-homomorphism $\rho: R \to \mathcal{S}$
given by $r \mapsto \rho_r$, right multiplication by $r$. These are sometimes called
the source map $\lambda$ and the target map $\rho$ of the bialgebroid. These two mappings commute within $\mathcal{S}$ at all values:
$$ \rho_r \circ \lambda_s = \lambda_s \circ \rho_r, \ \ r,s \in R := C_A(B),
$$
whence we may define a bimodule by composing strictly from the left: $$ {}_R\mathcal{S}_R:\ \ r \cdot \alpha \cdot s = \lambda_r \circ \rho_s \circ \alpha = r\alpha(?)s, \ \ \alpha \in \mathcal{S},\ r,s \in R.$$

Next we equip $\mathcal{S}$ with an $R$-coring structure $(\mathcal{S},R,\cop, \eps)$ as follows.  The comultiplication
$\cop: \mathcal{S} \to \mathcal{S} \o_R \mathcal{S}$ is an $R$-$R$-homomorphism given by
\begin{equation}
\label{eq: so}
\cop(\alpha) := \sum_i \alpha( ? t_i^1)t_i^2 \o_R \beta_i = \sum_j \gamma_j \o_R u_j^1 \alpha(u_j^2 ?)
\end{equation}
in terms of left D2 quasibases in the first equation or right D2 quasibases in the second equation.  
There is a simplification that shows this comultiplication is a generalization of the one in \cite[Lu]{Lu} (for the linear
endomorphisms of an algebra):
\begin{equation}
\label{eq: tys9}
\phi:\, \mathcal{S} \o_R \mathcal{S} \stackrel{\cong}{\longrightarrow} \Hom ({}_BA \o_B A_B, {}_BA_B)
\ \ \ \phi(\alpha \o_R \beta)(x \o_B y) = \alpha(x)\beta(y)
\end{equation}
 for all $\alpha, \beta \in \mathcal{S}$ and $x,y \in A$.
From a variant of eq.~(\ref{eq: tys}) we obtain:
\begin{equation}
\label{eq: lu}
 \phi(\cop(\alpha))(x \o y) =  \alpha(xy),\ \ \ (x,y \in A,\ \alpha \in \End {}_BA_B). 
\end{equation}
The counit $\eps: \mathcal{S} \to R$ is given by evaluation at the unity element,
$\eps(\alpha) = \alpha(1_A)$, again an $R$-$R$-homomorphism.  It is then apparent from eq.~(\ref{eq: lu}) that
$$(\eps \o_R \id_{\mathcal{S}})\cop(\alpha) = \eps(\alpha\1)\cdot \alpha\2 = \lambda_{\alpha\1(1)}\alpha\2 = \alpha, \ \
(\alpha \in \mathcal{S})$$ 
using a reduced Sweedler notation for the coproduct of an element,
and a similar equation corresponding to $(\id_{\mathcal{S}} \o_R \eps)\cop = \id_{\mathcal{S}}$.  

Finally, the comultiplication and counit satisfy additional bialgebra-like axioms that make $(\mathcal{S},R,\lambda, \rho,\cop, \eps)$
a left $R$-bialgebroid \cite[p.\ 80]{KS}.  These are:
\begin{equation}
\eps(1_{\mathcal{S}}) = 1_R,
\end{equation}
which is obvious,
\begin{equation}
\cop(1_{\mathcal{S}}) = 1_{\mathcal{S}} \o_R 1_{\mathcal{S}},
\end{equation}
which follows from eq.~(\ref{eq: lu}), 
\begin{equation}
\label{eq: pre}
\forall \, \alpha \in \mathcal{S},\ \alpha\1 \circ \rho_r \o_R \alpha\2 = \alpha\1 \o_R \alpha\2 \circ \lambda_r
\end{equation}
which follows from the equation defining $\phi$ (since both sides yield $\alpha(xry)$),
\begin{equation}
\cop(\alpha \circ \beta) = \cop(\alpha) \cop(\beta) = \alpha\1 \circ \beta\1 \o_R \alpha\2 \circ \beta\2 
\end{equation}
where eq.~(\ref{eq: pre}) justifies the use of a tensor algebra product in $\Im \cop \subseteq \mathcal{S} \o_R \mathcal{S}$
and the equation follows again from equation defining $\phi$ (as both sides equal $\alpha(\beta(xy))$), and at last
the easy 
\begin{equation}
\eps(\alpha \circ \beta) = \eps(\alpha \circ \lambda_{\eps(\beta)}) = \eps(\alpha \circ \rho_{\eps(\beta)}).
\end{equation}

On occasion the bialgebroid $\mathcal{S}$ is a \textit{Hopf algebroid} \cite{BS}, i.e. possesses an \textit{antipode} $\tau: \mathcal{S} \to \mathcal{S}$.  This
is \textit{anti-automorphism} of the algebra $\mathcal{S}$ which satisfies:
\begin{enumerate}
\item \label{eq: HA1} $\tau \circ \rho = \lambda$,
\item \label{eq: HA2} $\tau^{-1}(\alpha\2)\1 \o_R \tau^{-1}(\alpha\2)\2 \alpha\1 = \tau^{-1}(\alpha) \o 1_{\mathcal{S}}$ \ \ ($\forall\, \alpha \in \mathcal{S}$);
\item \label{eq: HA3} $\tau(\alpha\1)\1 \alpha\2 \o_R \tau(\alpha\1)\2 = 1_{\mathcal{S}} \o_R \tau(\alpha)$ \ \ ($\forall\, \alpha \in \mathcal{S}$).
\end{enumerate}

Examples of Hopf algebroids are weak Hopf algebras \cite{ENO} 
and Hopf algebras, including group algebras and
enveloping algebras of Lie algebras \cite{Mo}.  Lu \cite{Lu}
provides the example $A \o A^{\rm op}$ of a Hopf algebroid over any algebra $A$ with twist being the antipode, and another bialgebroid the linear endomorphisms $\End A$,  which is a particular case of 
the construction $\mathcal{S}$.  

A homomorphism of $R$-bialgebroids $\mathcal{S}' \to \mathcal{S}$ 
is an algebra homomorphism which commutes with the source and target mappings,  which
additionally is an $R$-coring homomorphism (so there are three commutative triangles and a commutative square
for such an algebra homomorphism to satisfy) \cite{LK2003}.   If $\mathcal{S}'$ and $\mathcal{S}$ additionally
come equipped with antipodes $\tau'$ and $\tau$, respectively, then the homomorphism
$\mathcal{S}' \to \mathcal{S}$ is additionally a homomorphism of Hopf algebroids if it commutes with the antipodes
in an obvious square diagram.

\subsection{Pseudo-Galois extensions.}
If $\sigma$ is automorphism of the algebra $A$, we let
$A_{\sigma}$ denote the bimodule $A$ twisted on the right by
$\sigma$, with module actions
defined by $x \cdot a \cdot y = xa \sigma(y)$ for $x,y, a \in A$.
Two such bimodules $A_{\sigma}$ and $A_{\tau}$ twisted
by automorphisms $\sigma, \tau: A \to A$ are $A$-$A$-bimodule isomorphic if and only if there is an invertible
element $u \in A$ such that $\tau \circ \sigma^{-1}$ is
the inner automorphism by $u$ (for send $1 \mapsto u$).

Let  $B$ be a subalgebra of $A$. Recall the 
characterization of a group-Galois extension 
 $A \| B$, where only two conditions need be met.  First, 
there is a finite group $G$ of automorphisms of $A$ such
that $B = A^G$, i.e.\ the elements of $B$ are fixed under each
automorphism of $G$ and each element of $A$ in the complement of $B$ is moved by some automorphism of $G$.  Second, there are
elements $a_i, b_i \in A$, $i = 1,\ldots,n$ such
that $\sum_i a_i b_i = 1$ and $\sum_i a_i \sigma(b_i) = 0$ if $\sigma \not = \id_A$.  
 
Since $E: A \to B$ defined by $E(a) = \sum_{\sigma \in G} \sigma(a)$
is a Frobenius homomorphism with dual bases $a_i, b_i$,
it follows that there is an $A$-$A$-bimodule isomorphism
between the tensor-square and the semi-direct product of $A$ and $G$:
\begin{equation}
\label{eq: ach}
h: \ A \o_B A \longrightarrow A \rtimes G, \ \ h(x \o y) := \sum_{\sigma \in G}^{} x \sigma y = \sum_{\sigma \in G} x \sigma(y) \sigma 
\end{equation}
(For the inverse is given by $h^{-1}(a \tau) = \sum_i a \tau(a_i) \o_B b_i$.)
Thus $A \o_B A$ is isomorphic to $\oplus_{\sigma \in G} A_{\sigma}$
as $A$-$A$-bimodules. Mewborn and McMahon \cite{MM} relax this condition as follows:

\begin{definition}
The algebra $A$ is a \textit{pseudo-Galois} extension of a subalgebra
$B$ if there is a finite set $\mathcal{G}$ of $B$-automorphisms (i.e.
fixing elements of $B$) and a positive integer $N$ such that
$A \o_B A$ is isomorphic to a direct summand of $\oplus_{\sigma \in \mathcal{G}} A_{\sigma}^N$: in symbols this becomes
\begin{equation}
\label{eq: pseudo-Galois}
A \o_B A \, \oplus \, * \, \cong \, \oplus_{\sigma \in \mathcal{G}} A_{\sigma}^N, 
\end{equation}
in terms of $A$-$A$-bimodules. Assume with no loss of generality that $\mathcal{G}$ is minimal
in the group of $B$-automorphisms $\Aut_B(A)$ with respect to this property and $A_{\sigma} \not \cong A_{\tau}$
if $\sigma \not = \tau$ in $\mathcal{G}$.  
\end{definition}

It is clear that Galois extensions are pseudo-Galois.  For example, let $A$ be a simple  ring with finite group $G$ of outer automorphisms of $A$, then $A$ is Galois over
its fixed subring $B = A^G$, cf.\ \cite[2.4]{Mo80}. 
Another example of a pseudo-Galois extension is an \textit{H-separable
extension} $A \| B$, which by definition satisfies $A \o_B A \oplus * \cong A^N$ as
$A$-$A$-bimodules, so we let $\mathcal{G} = \{ \id_A \}$ in the definition
above \cite{NEFE, LK2003, MM}. For example, if $A$ is a simple ring, $G$ a finite group
of outer automorphisms of $A$ such that each nonidentity automorphism moves an element
of the center of $A$, then the skew group ring $A \rtimes G$ is H-separable over $A$ \cite{MS}. 

Note that the definition of pseudo-Galois extension $A \| B$ leaves open the possibility
that $B$ is a proper subset of the invariant subalgebra $A^{\mathcal{G}}$: if $B \subset C \subset
A^{\mathcal{G}}$ and $C$ is a separable extension of $B$, then $A \| C$ is also a pseudo-Galois
extension, since one may show that $A \o_C A \oplus * \cong A \o_B A$ as natural $A$-$A$-bimodules via  the separability
element. 
 Conversely, if $A \| C$ is a pseudo-Galois extension and $C \supset B$ is H-separable,
then $A \| B$ is pseudo-Galois by noting that $A \o_B A \cong$ $ A \o_C C \o_B C \o_C A$. 
For example, if $E \| F$ is a finite Galois extension of fields where $F$ is the quotient
field of a domain $R$, then the ring extension $E \| R$ is pseudo-Galois.

In the next proposition, we note that pseudo-Galois extensions are depth two by means of a  characterization of pseudo-Galois extensions using
\textit{pseudo-Galois elements}.

\begin{prop}
 An algebra extension $A \| B$ is pseudo-Galois iff
there is a finite set $\mathcal{G}$ of $B$-automorphisms 
and $N$ elements $r_{i, \sigma} \in C_A(B)$
and $N$ elements $e_{i,\sigma} \in ({}_{\sigma}A \o_B A)^A$ for each element $\sigma \in \mathcal{G}$ satisfying
\begin{equation}
\label{eq: pseudo-Galois elements}
1 \o_B 1 = \sum_{\sigma \in \mathcal{G}} \sum_{i=1}^N r_{i, \sigma}
e_{i, \sigma}
\end{equation}
As a consequence, $A \| B$ is left and right D2 with left
and right D2 quasibases derived from the elements $r_{i, \sigma}$
and $e_{i, \sigma} \in (A \o_B A)^B$.  
\end{prop}

\begin{proof}
($\Rightarrow$) We note that the condition~(\ref{eq: pseudo-Galois})
implies the existence of $N$ pairs of mappings $f_{i, \sigma}$ and 
$g_{i,\sigma}$ for each $B$-automorphism $\sigma \in \mathcal{G}$ satisfying 
\begin{equation}
\label{eq: hirata}
\sum_{i=1}^N \sum_{\sigma \in \mathcal{G}}
g_{i, \sigma} \circ f_{i, \sigma} = \id_{A \o_B A}.
\end{equation}
  The mappings
simplify as 
$$ f_{i, \sigma} \in \Hom ({}_AA \o_B A_A, {}_A(A_{\sigma})_A) \cong C_A(B),  $$
via $F \mapsto F(1 \o 1)$ with inverse  $r \mapsto (a \o a' \mapsto
a r \sigma(a'))$, as well as mappings
$$ g_{i, \sigma} \in \Hom ({}_A(A_{\sigma})_A, {}_AA \o_B A_A) \cong
({}_{\sigma}A \o_B A)^A, $$
via $f \mapsto f(1)$ with inverse $e \mapsto (a \mapsto ae)$
where $e \in ({}_{\sigma}A \o_B A)^A$ iff $ea = \sigma(a)e$
for each $a \in A$.  

If $r_{i, \sigma}$ corresponds via the isomorphism above with
$f_{i, \sigma}$, then $f_{i, \sigma}(x \o_B y) = x r_{i, \sigma} \sigma(y)$ for each $x,y \in A$.
If $e_{i, \sigma}$ corresponds via the other isomorphism above
with $g_{i, \sigma}$, then $g_{i, \sigma}(a) = a e_{i, \sigma}$.
We compute:
\begin{equation}
\label{eq: needed}
 x \o_B y = \sum_{i, \sigma \in \mathcal{G}}
(g_{i, \sigma} \circ f_{i, \sigma})(x \o y) = \sum_{\sigma \in \mathcal{G}} \sum_{i = 1}^N xr_{i, \sigma} \sigma(y) e_{i, \sigma},
\end{equation}
which shows that 
\begin{equation}
\label{eq: right D2 pseudo basis}
\lambda_{r_{i, \sigma}} \circ \sigma \in \End {}_BA_B,  \ \ 
e_{i, \sigma} \in (A \o_B A)^B
\end{equation}
 are right D2 quasibases.  
Setting $x = y = 1$ we obtain the eq.~(\ref{eq: pseudo-Galois elements}). Finally use the twisted centralizer property
$ea = \sigma(a) e$ for $a \in A$ and $e \in ({}_{\sigma}A \o_B A)^A$
to obtain
$$ x \o_B y = \sum_{\sigma \in \mathcal{G}} \sum_{i=1}^N e_{i, \sigma}
\sigma^{-1}(x) \sigma^{-1}(r_{i, \sigma})y. $$
Hence, the following are left D2 quasibases for $A \| B$:
\begin{equation}
\label{eq: left D2 pseudo basis}
e_{i, \sigma} \in (A \o_B A)^B, \ \ \sigma^{-1} \circ \rho_{r_{i, \sigma}} \in \End {}_BA_B.  
\end{equation}

($\Leftarrow$) Conversely, suppose we are given a finite set $\mathcal{G}$ of
$B$-automorphisms, and for each $\sigma \in \mathcal{G}$, $N$ centralizer elements $r_{i, \sigma} \in C_A(B)$,
and $N$ twisted $A$-central elements $e_{i, \sigma} \in ({}_{\sigma}A \o_B A)^A$
for $i = 1, \ldots,N$ such that eq.~(\ref{eq: pseudo-Galois elements}) holds.
By multiplying the equation from the left by $x \in A$ and from the right by $y \in A$,
we obtain eq.~(\ref{eq: needed}) and then eq.~(\ref{eq: hirata}) by defining $f_{i, \sigma}$
and $g_{i,\sigma}$ as before, which is of course equivalent to the condition~(\ref{eq: pseudo-Galois})
for pseudo-Galois extension. 
\end{proof}

We note that pseudo-Galois elements in eq.~(\ref{eq: pseudo-Galois elements}) specialize to H-separability elements in
case $\mathcal{G} = \{ \id_A \}$ \cite[2.5]{NEFE}. 

Recall that in homological algebra the enveloping algebra of an algebra
$A$ is denoted by $A^e := A \o A^{\rm op}$.
 
\begin{theorem}
\label{th-epic}
Suppose $A \| B$ is pseudo-Galois extension satisfying condition~(\ref{eq: pseudo-Galois})
with $G$ the subgroup generated by $\mathcal{G}$ within $\Aut_B A$ and $R$ the centralizer $C_A(B)$.
Then there is a Hopf algebroid denoted by $R^e \bowtie G$ which maps epimorphically as $R$-bialgebroids
onto the left bialgebroid $\mathcal{S} = \End {}_BA_B$.  
\end{theorem}

\begin{proof}
Denote the identity in $G$ by $e$ and the canonical anti-isomorphism $R \to R^{\rm op}$ by
$r \mapsto \overline{r}$ satisfying $\overline{r}\, \overline{s} = \overline{sr}$
for every $r,s \in R$.  Note that the $B$-automorphisms of $G$ restrict to automorphisms
of the centralizer $R$. The notation $C$ for $R^e \bowtie G$ is adopted at times.   
The Hopf algebroid structure of $R^e \bowtie G$ is given by
\begin{enumerate}
\item as a $K$-module (over the ground ring $K$) $R^e \bowtie G = R \o R^{\rm op} \o K[G]$
where $K[G]$ is the group $K$-algebra of $G$;
\item multiplication  given by 
\begin{equation}
\label{eq: prod}
(r \o \overline{s} \b \sigma)(u \o \overline{v} \b \tau) = r \sigma(u) \o \overline{v \tau^{-1}(s)} \b \sigma \tau 
\end{equation}
with unity element $1_C = 1_R \o \overline{1_R} \b e$,
\item source map $s_L: R \to R^e \bowtie G $ given by $s_L(r) = r \o \overline{1} \bowtie e$,
\item target map $t_L: R^{\rm op} \to R^e \bowtie G$ given by $t_L(r) = 1 \o \overline{r} \b e$,
\item counit $\eps_C: R^e \bowtie G \to R$ given by $\eps_C(r \o \overline{s} \b \sigma) = r \sigma(s)$,
\item comultiplication $\cop_C:  C \to C\o_R C$ is given by 
\begin{equation}
\label{eq: coprod}
\cop(r \o \overline{s} \b \sigma) =
(r \o \overline{1} \b \sigma ) \o_R (1 \o \overline{s}\b \sigma ),
\end{equation}
\item and the antipode $\tau: C \to C$ by $\tau(r \o \overline{s} \b \sigma) = s \o \overline{r} \b \sigma^{-1}$
\end{enumerate}
We will postpone the proof that this defines a Hopf algebroid over $R$ until the next section
where it is shown  more generally for an involutive Hopf algebra $H$ and its 
$H$-module algebras.  

The epimorphism of left $R$-bialgebroids $\Psi: R^e \b G \to \mathcal{S}$ is given by
\begin{equation}
\label{eq: epi}
\Psi(r \o \overline{s} \b \sigma) = \lambda_r \circ \sigma \circ \rho_s
\end{equation}
We note that $\Psi$ is an algebra homomorphism by comparing eq.~(\ref{eq: prod}) with
$$ \lambda_r \circ \sigma \circ \rho_s \circ \lambda_u \circ \tau \circ \rho_v =
\lambda_{r \sigma(u)} \circ \sigma \circ \tau \circ \rho_{v \tau^{-1}(s)}, $$
and $\Psi(1_C) = 1_{\mathcal{S}}$.  
The mapping $\Psi$ is epimorphic since each $\beta \in \mathcal{S}$ may be expressed as
a sum of mappings of the form $\lambda_r \circ \sigma \circ \rho_s$
where $\sigma \in G$ and $r,s \in R$.  To see this, apply
$\mu (\id_A \o \beta)$ to eq.~(\ref{eq: needed}) with $x = 1$,
which yields
$$ \beta(y) = \sum_{i, \sigma \in G} r_{i, \sigma} \sigma(y) e_{i, \sigma}^1 \beta(e_{i, \sigma}^2)$$
where $e_{i, \sigma}^1 \beta(e_{i, \sigma}^2) \in R$ for each $i$ and $\sigma$.  

Note next that $\Psi$ commutes with source, target and counit maps.
For $\Psi(s_L(r)) = \Psi(r \o \overline{1} \b e) = \lambda_r$ and 
$\Psi(t_L(r)) = \Psi(1 \o \overline{r} \b e) = \rho_r$ for $r \in R$
(so  $\Psi: C \to \mathcal{S}$ is an $R$-$R$-bimodule map).
The map $\Psi$ is counital since $$ \eps(\Psi(r \o \overline{s} \b \sigma)) = \lambda_r (\sigma(\rho_s(1))) = r \sigma(s) =
\eps_C(r \o \overline{s} \b \sigma). $$

Using the isomorphism $\phi: \mathcal{S} \o_R \mathcal{S} \to \Hom ({}_BA \o_B A_B, {}_BA_B)$ for a depth two 
extension $A \| B$ 
defined as above by $\phi(\alpha \o_R \beta)(x \o_B y) = \alpha(x)\beta(y)$,
note from eq.~(\ref{eq: coprod}) that $\Psi$ is comultiplicative:  
$$\phi( (\Psi \o_R \Psi) (\cop_C(r \o \overline{s} \b \sigma)))(x \o_B y) = 
 \lambda_r(\sigma(x)) \sigma(\rho_s(y)) = $$
$$ \phi(\cop(\lambda_r \circ \sigma \circ \rho_s))(x \o_B y) = \phi (\cop(\Psi(r \o \overline{s} \b \sigma)))(x \o_B y)
$$ since $\sigma \in G$ is a group-like element satisfying $\sigma\1 \o_R \sigma\2 = \sigma \o_R \sigma$ (corresponding to the automorphism condition).  
\end{proof}
 
\begin{cor}
If the algebra extension $A \| B$ is H-separable, then $G = \{ \id_A \}$
and $\Psi: R^e \to \mathcal{S}$ is
an isomorphism of bialgebroids, whence $\mathcal{S}$ has the antipode $\Psi \circ \tau \circ \Psi^{-1}$.
  If $A \| B$ is $G$-Galois, then $\Psi: R^e \b G \to \mathcal{S}$ is a split epimorphism of bialgebroids.
\end{cor}

\begin{proof}
Note that $R^e$ is isomorphic as algebras to the subalgebra $R^e \b \{e \}$.  
The first statement follows from \cite{LK2003}, since $\Psi(r \o \overline{s}) = \lambda_r \circ
\rho_s$ is shown there to be an isomorphism of bialgebroids.

 If $A \| B$ is $G$-Galois, then $A \o_B A \cong
A \rtimes G$ via $h$ above.  Since each $\sigma \in G$ fixes elements of $B$, it follows
that $(A \o_B A)^B \cong R \rtimes G$. 
Since $A \| B$ is a Frobenius extension, $\End A_B \cong A \o_B A$ via
$f \mapsto \sum_i f(a_i) \o b_i$ with inverse $x \o_B y \mapsto \lambda_x \circ E \circ \lambda_y$.
This restricts to $\End {}_BA_B \cong (A \o_B A)^B$. Putting the two together yields
$$\Phi:\ \mathcal{S} \stackrel{\cong}{\longrightarrow} R \rtimes G, \ \ \ \Phi(\alpha) := \sum_{\sigma \in G} \sum_{i=1}^n \alpha(a_i) \sigma(b_i)\sigma$$ with inverse $r \rtimes \tau \mapsto \lambda_r \circ \tau$. Then the algebra epimorphism $\Phi \circ \Psi:\ R^e \b G \to R \rtimes G$ simplifies to $$ (\Phi \circ \Psi)(r \o \overline{s} \b \sigma) = \Phi(\lambda_r \circ \sigma \circ \rho_s) = \sum_{\tau \in G} \sum_{i=1}^n \lambda_{r\sigma(a_i s \tau(b_i))} \circ \sigma \circ \tau.
$$  which is split by the monomorphism $R \rtimes G \to R^e \b G$ given
by $r \rtimes \sigma \mapsto r \o \overline{1} \b \sigma$, an algebra
homomorphism by an application of eq.~(\ref{eq: prod}).   
\end{proof}


\section{An enveloping Hopf algebroid over algebras in  certain tensor categories}
 
Let $H$ be a Hopf algebra with bijective antipode $S$ and $A$ a left $H$-module
algebra, i.e. an algebra in the tensor category of $H$-modules.   
Motivated by the left bialgebroid of a pseudo-Galois extension
as studied in section~3, we define a
type of enveloping algebra $A^e \b H$ for the smash product algebra $A \rtimes H$.  It is
a left bialgebroid over $A$, and a Hopf algebroid in case $H$ is involutive such as
a group algebra or the enveloping algebra of Lie algebra.  In terms of noncommutative algebra, it is the minimal algebra
which contains subalgebras isomorphic to the
Hopf algebra $H$, the standard enveloping algebra $A^e$
of an algebra $A$, 
and the 
semi-direct or crossproduct algebra $A \rtimes H$ as
well as its derived right crossproduct algebra $H \ltimes A^{\rm op}$. 
In terms of category theory, it is derived from the pushout construction \cite{MAC} of the inclusion
$H \into A \rtimes H$ and its opposite via the isomorphism $S: H \to H^{\rm cop,\, op}$.  

\begin{theorem}
\label{th-bowtie}
  Suppose $B := A^e \bowtie H$ is the
vector space $A \o A^{\rm op} \o H$ with multiplication
\begin{equation}
\label{eq: product}
(a \o \overline{b} \b h)( c \o \overline{d} \b k) :=
a (h\1 \cdot c) \o \overline{d(S(k\2) \cdot b)} \b h\2 k\1.
\end{equation}
Then $B$ is a left bialgebroid over $A$ with
structure given in eqs.~(\ref{eq: un}) through~(\ref{eq: quatre}).  If $S^2 = \id_H$, then
$B$ is a Hopf algebroid with antipode eq.~(\ref{eq: six}).  
\end{theorem}
\begin{proof}

Clearly the unity element $1_B = 1_A \o \overline{1_A} \b 1_H$.  
The multiplication is associative, since
$$
[(a \o \overline{b} \b h)( c \o \overline{d} \b k)]
(e \o \overline{f} \b j)  =  (a (h\1 \cdot c) \o \overline{d(S(k\2) \cdot b)} \b h\2 k\1 )(e \o \overline{f} \b j) = $$
$$  a(h\1 \cdot c)(h\2 k\1 \cdot e) \o \overline{f(S(j\3)\cdot d)(S(k\3 j\2)\cdot b)} \b h\3 k\2 j\1 = $$
$$
(a \o \overline{b} \b h)(c(k\1 \cdot e) \o \overline{f(S(j\2) \cdot d)} \b k\2 j\1  =  
(a \o \overline{b} \b h)[( c \o \overline{d} \b k)
(e \o \overline{f} \b j) ].
$$
It follows that $B$ is an algebra.

Define a source map $s_L: A \to B$ and target
map $t_L: A \to B$ by
\begin{equation}
\label{eq: un}
s_L(a) = a \o \overline{1_A} \b 1_H
\end{equation}
\begin{equation}
\label{eq: deux}
t_L(a) = 1_A \o \overline{a} \b 1_H,
\end{equation}
an algebra homomorphism and anti-homomorphism, respectively. It is evident that $t_L(x)s_L(y) = s_L(y)
t_L(x)$ for all $x,y  \in A$.  
The $A$-$A$-bimodule structure induced from
$x \cdot b \cdot y = s_L(x) t_L(y) b$ for $b \in B$
is then given by ($a,c \in A$, $h \in H$)
\begin{equation}
\label{eq: bi}
x \cdot (a \o \overline{c} \b h) \cdot y =
xa \o \overline{c(S(h\2)\cdot y)} \b h\1.
\end{equation}

The counit $\eps: B \to A$ is defined by
\begin{equation}
\label{eq: trois}
\eps(a \o \overline{c} \b h) := a (h \cdot c).
\end{equation}
Note that $\eps$ is an $A$-$A$-bimodule homomorphism
via its application to the RHS of eq.~(\ref{eq: bi}):
$$ \eps(xa \o \overline{c(S(h\2)\cdot y)} \b h\1)
= xa h\1 \cdot (c(S(h\2)\cdot y)) = xa(h \cdot c) y =
x \eps(a \o \overline{c} \b h) y, $$ 
since $h \cdot (xy) = (h\1 \cdot x)(h\2 \cdot y)$, 
the measuring axiom on $A$.  

The comultiplication $\cop: B \to B \o_A B$ is defined by
\begin{equation}
\label{eq: quatre}
\cop(a \o \overline{c} \b h) := (a \o \overline{1_A} \b h\1) \o_A (1_A \o \overline{c} \b h\2 ).
\end{equation}
It is an $A$-$A$-homomorphism:
$$ \cop(xa \o \overline{c(S(h\2)\cdot y)} \b h\1) =
(xa \o \overline{1_A} \b h\1) \o_A (1_A \o \overline{c(S(h\3) \cdot y)} \b h\2  $$
$$= x \cdot \cop(a \o \overline{c} \b h) \cdot y$$
by eq.~(\ref{eq: bi}). The left counit equation 
$(\eps \o_A \id_B)\cop = \id_B$ follows from
$$ \eps(a \o \overline{1} \b h\1) \cdot (1 \o \overline{c} \b h\2) = $$
$$ a(h\1 \cdot 1) \cdot (1 \o \overline{c} \b h\2) = a \o\overline{c} \b h, $$
since $h \cdot 1_A = \eps(h) 1_A$ in the $H$-module
algebra $A$.  The right counit equation
$(\id_B \o \eps)\cop = \id_B$ follows from
$$(a \o \overline{1} \b h\1) \cdot \eps(1 \o \overline{c} \b h\2) = (a \o \overline{1} \b h\1) \cdot (h\2 \cdot c) = $$
$$ a \o \overline{(S(h\2)h\3 \cdot c)} \b h\1 = a \o \overline{c} \b h. $$
Hence $(B,A,\cop, \eps)$ is an $A$-coring.

We check the remaining bialgebroid axioms:
$$ \cop(1_B) = 1_B \o 1_B, \ \ \ \eps(1_B) = 1_A  $$
are apparent from eqs.~(\ref{eq: quatre}) and~(\ref{eq: trois}).
The axiom corresponding to eq.~(\ref{eq: pre}) computes as:
$$ \cop(a \o b \b h)(t_L(c) \o_A 1_B) = (a \o \overline{1} \b h\1)(1 \o \overline{c} \b 1_H)
\o_A (1 \o \overline{b} \b h\2) $$
$$ =  (a \o \overline{c} \b h\1) \o_A (1 \o \overline{b} \b h\2), $$
 and on the other hand
$$ \cop(a \o \overline{b} \b h)(1_B \o_A s_L(c)) = (a \o \overline{1} \b h\1) \o_A (1 \o \overline{b} \b h\2)(c \o \overline{1} \b 1_H) = $$
$$ (a \o \overline{1} \b h\1) \o_A ( h\2 \cdot c \o \overline{b} \b h\3) =
(a \o \overline{S(h\2)h\3 \cdot c} \b h\1) \o_A (1 \o \overline{b} \b h\4) $$
which equals the RHS expression for $\cop(a \o b \b h)(t_L(c) \o_A 1_B)$.

Next, the comultiplication is multiplicative:
$$\cop((a \o \overline{b} \b h)(c \o \overline{d} \b k)) = \cop(a (h\1 \cdot c) \o \overline{d(S(k\2) \cdot b)} \b h\2 k\1) = $$
$$ (a (h\1 \cdot c) \o \overline{1} \b h\2 k\1)
 \o_A (1_A \o \overline{d( S(k\3) \cdot b)} \b h\3k\2 = $$
$$ (a \o \overline{1} \b h\1)(c \o \overline{1} \b k\1) \o_A (1 \o \overline{b} \b h\2)(1 \o \overline{d} \b k\2) 
 = \cop(a \o \overline{b} \b h) \cop(c \o \overline{d} \b k). $$
The counit satisfies 
$$ \eps((a \o \overline{b} \b h)(c \o \overline{d} \b k)) = \eps(a (h\1 \cdot c) \o \overline{d(S(k\2) \cdot b)} \b h\2 k\1) = $$
$$ a (h\1 \cdot c)(h\2 k\1 \cdot (d (S(k\2) \cdot b))) = a(h\1 \cdot c) (h\2 k \cdot d)(h\3 \cdot b) = $$
$$ = \eps((a \o \overline{b} \b h)(c (k \cdot d) \o \overline{1} \b 1_H)) = \eps((a \o \overline{b} \b h) s_L(\eps(c \o \overline{d} \b k))).$$
Similarly, $\eps((a \o \overline{b} \b h)(c \o \overline{d} \b k)) = 
\eps((a \o \overline{b} \b h) t_L(\eps(c \o \overline{d} \b k)))$ for
all $a,b,c,d \in A$, $h,k \in H$.  Thus $B$ is a bialgebroid over $A$.

Suppose the antipode on $H$ is bijective and satisfies $S^2 = \id_H$. 
Define an antipode on $B$ by ($a,b \in A$, $h \in H$) 
\begin{equation}
\label{eq: six}
\tau(a \o \overline{b} \b h) = b \o \overline{a} \b S(h)
\end{equation}
Denote the compositional inverse of $S$ by $\overline{S}$.  
Then $\tau$ has inverse,
$$ \tau^{-1}(a \o \overline{b} \b h) = b \o \overline{a} \b \overline{S}(h). $$

Note that $\tau$ is an anti-automorphism of $B$:
$$ \tau(c \o \overline{d} \b k)\tau(a \o \overline{b} \b h) = (d \o \overline{c} \b S(k))(b \o \overline{a} \b S(h)) = $$
$$ d (S(k\2) \cdot b) \o \overline{a(S^2(h\1)\cdot c)} \b S(k\1)S(h\2) =
\tau(a (h\1 \cdot c) \o \overline{d(S(k\2) \cdot b)} \b h\2 k\1) $$
$$ = \tau((a \o \overline{b} \b h)(c \o \overline{d} \b k)). $$

The antipode satisfies the three axioms~(\ref{eq: HA1})-(\ref{eq: HA3}):
$$ \tau(t_L(a)) = \tau(1 \o \overline{a} \b 1_H) = a \o \overline{1} \b 1_H = s_L(a), $$
for all $a \in A$.  
Next, for $b := a \o \overline{c} \b h \in B$, 
$$ \tau^{-1}(b\2)\1 \o_A \tau^{-1}(b\2)\2 b\1 =$$
$$ \tau^{-1}(1 \o \overline{c} \b h\2)\1 \o_A \tau^{-1}(1 \o \overline{c} \b h\2)\2 (a \o \overline{1} \b h\1) = $$
$$(c \o \overline{1} \b \overline{S}(h\4)) \o_A
(\overline{S}(h\3) \cdot a \o \overline{1} \b \overline{S}(h\2)h\1 = $$
$$ c \o \overline{S(\overline{S}(h\2)\overline{S}(h\1)\cdot a}
\b \overline{S}(h\3) \o_A 1_B =  (c \o \overline{a} \b \overline{S}(h)) \o_A 1_B = \tau^{-1}(b) \o_A 1_B .$$
Continuing our notation $b = a \o \overline{c} \b h \in B$, note too that
$$ \tau(b\1)\1 b\2 \o_A \tau(b\1)\2 = (1 \o \overline{1} \b S(h\2))(1 \o \overline{c} \b h\3) \o_A (1 \o \overline{a} \b S(h\1))  $$
$$ = (1 \o \overline{c} \b 1_H) \o_A (1 \o \overline{a} \b S(h)) = 1_B \o_A (c \o \overline{a} \b S(h)) = 1_B \o_A \tau(b). $$
Hence, $B$ is a Hopf algebroid.
\end{proof}

Given a group $G$, its group algebra $K[G]$ over a commutative ring $K$ is
an involutive Hopf algebra \cite{Mo}.  Moreover, if $G$ acts by automorphisms
on a $K$-algebra $A$, then $A$ is a left $K[G]$-module algebra and $A \rtimes G$
is identical with the semidirect product \cite{Mo}.  Thus the construction $R^e \b G$ (covering
the left bialgebroid of a pseudo-Galois extension in Section~3) is a Hopf algebroid,
and we record the following.
\begin{cor}
Given a $K$-algebra $A$ and a group $G$ of algebra automorphisms of $A$,
the algebra $A^e \b K[G]$ is a Hopf algebroid over $A$.
\end{cor}

Recall that Lu \cite{Lu} defines over an algebra $A$ a Hopf algebroid $A^e$.  This
is a Hopf subalgebroid of the construction in the theorem above.
\begin{cor}
\label{cor-jay}
Let $H$ be an involutive Hopf algebra and $A$ a left $H$-module algebra.
Then the Hopf algebroid $A^e \b H$ contains subalgebras isomorphic to
\begin{enumerate}
\item Lu's Hopf algebroid $A\o A^{\rm op}$
\item the semidirect product $A \rtimes H$
\item its derived right crossproduct $H \ltimes A^{\rm op}$
\item the Hopf algebra $H$
\end{enumerate}
\end{cor}
\begin{proof}
It is easy to check from eq.~(\ref{eq: product}) that the following mappings
\begin{enumerate}
\item $ A^e \into A^e \b H$ given by $a \o \overline{b} \mapsto a \o \overline{b} \b 1_H$
is an algebra monomorphism as well as a homomorphism of Hopf algebroids over $A$
(i.e., it commutes with the source, target, counit, comultiplication and antipode maps
above and those given in \cite{Lu});
\item $j_1: A \rtimes H \into A^e \b H$ given by $j_1(a \# h) := a \o \overline{1_A} \b h$
is an algebra monomorphism, where we recall
 that the multiplication in $A \rtimes H$ is given
by $$(a \# h)(b \# k) = a(h\1 \cdot b) \# h\2 k $$
\item $j_2: H \ltimes A^{\rm op} \into A^e \b H$ given by $j_2(h \# \overline{a}) := 
1 \o \overline{a} \b h$ is an algebra monomorphism, where $\overline{a} \cdot h :=
\overline{S(h) \cdot a}$ defines the derived right action of $H$ on $A^{\rm op}$
and the multiplication in $H \ltimes A^{\rm op}$
(cf.\ \cite[p.\ 22]{Maj}) is given by $$(h \# \overline{a})(k \# \overline{b}) =
hk\1 \#  (\overline{a} \cdot k\2) \overline{b}. $$
\item $H \into A^e \b H$ given by $h \mapsto 1 \o \overline{1} \b h$ is an algebra
monomorphism as well as a Hopf algebra homomorphism (for it commutes with the counit, comultiplication
and antipode mappings of $H$ and $A^e \b H$
if $A$ is  a faithful $K$-algebra and  $K$ is identified
with $K1_A$.).  
\end{enumerate}
\end{proof}

The construction $A^e \b H$ for a Hopf algebra and a left $H$-module algebra
is a type of cofibered sum \cite[p.\ 99]{Sp} of the algebra monomorphisms $\iota_1: H \into A \rtimes H$
and $\iota_2: H \into H \ltimes A^{\rm op}$  defined
by $\iota_1(h) := 1_A \# h$ and $\iota_2(h) := h \# \overline{1_A}$ for each $h \in H$.  
We note that $j_1 \circ \iota_1 = j_2 \circ \iota_2$, both sending $h \mapsto 1_A \o \overline{1_A} \b h$.  Also define the algebra monomorphism $k_1: A \into A \rtimes H$ by $k_1(a) :=
a \# 1_H$ and anti-monomorphism $k_2: A \into H \ltimes A^{\rm op}$ by
$k_2(a) = 1_H \# \overline{a}$. Note that $$ j_2(k_2(a)) j_1(k_1(b)) = (1_A \o \overline{a} \b 1_H) (b \o \overline{1_A} \b 1_H)= (b \o \overline{1_A} \b 1_H)(1_A \o \overline{a} \b 1_H) $$
$$ = j_1(k_1(b))j_2(k_2(a)), $$
for all $a, b \in A$.  

\begin{prop}
\label{prop-universal}
Suppose $B$ is an algebra with monomorphisms $f_1: A \rtimes H \into B$ and
$f_2: H \ltimes A^{\rm op} \into B$ such that $f_1 \circ \iota_1 = f_2 \circ \iota_2$ (i.e., satisfying the commutative square in Figure~1) 
and $f_1(k_1(a)) f_2(k_2(b)) = f_2(k_2(b))f_1(k_1(a))$ for all $a, b \in  A$.
Then there is a uniquely defined algebra homomorphism $F: A^e \b H \to B$ such that $F \circ j_i = f_i$ for $i = 1,2$.  
\end{prop}
\begin{figure}
$$\begin{diagram}
H &&  \rTo^{\iota_1}&& A \rtimes H \\
\dTo^{\iota_2} && && \dTo_{f_1}    \\
H \ltimes A^{\rm op} & & \rTo^{f_2}  && B
\end{diagram}$$
\caption{$A^e \b H$ is the cofibered sum of $\iota_{1,2}: H \into A \rtimes H,\, H \ltimes A^{\rm op}$ such that $A^e$ embeds homomorphically.}
\end{figure}
\begin{proof}
Define $F: A^e \b H \to B$ by 
\begin{equation}
F(a \o \overline{b} \b h) := f_1(a \# h) f_2(1_H \# \overline{b})
= f_1( a \# 1_H) f_2(h \# \overline{b}).
\end{equation}
  The second equality follows from $a \# h = (a \# 1_H)
(1_A \# h)$ and $f_1 \circ \iota_1 = f_2 \circ \iota_2$.  
It follows that $F \circ j_i = f_i$ for $i = 1,2$ since $f_i(1) = 1_B$.
Then the uniqueness of $F$ follows from noting $a \o \overline{b} \b h = 
(a \o \overline{1_A} \b h)(1 \o \overline{b} \b 1_H)$  and the homomorphic
property of $F$.  
 We compute that $F$ is an algebra homomorphism by using $f_1(k_1(c)) f_2(k_2(b)) = f_2(k_2(b))f_1(k_1(c))$ in the third equality:
$$ F(a \o \overline{b} \b h) F( c \o \overline{d} \b k) = f_1(a \# h)
f_2(1_H \# \overline{b}) f_1(c \# 1_H)f_2(k \# \overline{d}) = $$
$$ f_1(a \# h)f_1(c \# 1_H)f_2(1_H \# \overline{b})f_2(k \# \overline{d}) = f_1(a(h\1 \cdot c) \# h\2) f_2(k\1 \# \overline{d (S(k\2) \cdot b)} )  $$
$$= F((a \o \overline{b} \b h)(c \o \overline{d} \b k)), $$
by comparing the last equation with eq.~(\ref{eq: product}).  
\end{proof}

The homomorphism $F$ may fail to be monic as for example
when $A$ is a commutative algebra, $H$ acts trivially on $A$ and $B = A \o H$. 

\begin{example}
\begin{rm}
Enveloping algebras of Lie algebras are involutory Hopf algebras with comultiplication
defined via primitive elements and the antipode via sign change.  The Weyl algebra
$K[X,Z \, | \, XZ + 1 = ZX ]$ 
 is isomorphic to the semi-direct product $A \ltimes K[Z]$ of the one-dimensional Lie algebra
$K[Z]$ acting by Leibniz derivation on the one-variable polynomial algebra $A = K[X]$ \cite{Mo}.  
The enveloping Hopf algebroid is then the pushout of the inclusion 
$K[Z] \into K[X,Z \| XZ + 1 = ZX ]$ with itself:

\begin{equation}
K[X,Y] \bowtie K[Z] \cong K[X,Y,Z \|  XZ + 1 = ZX, \, XY = YX, \, YZ + 1 = ZY ]
\end{equation}
with Hopf algebroid comultiplication given on monomials by (integers $ k \geq p, q \geq 0$)
\begin{equation}
\cop(X^n Y^m Z^k) = \sum_{p+ q = k} \left( \begin{array}{c}
k \\
p
\end{array}
\right) X^n Z^p \o_A Y^m Z^q 
\end{equation}
counit by
\begin{equation}
\eps(X^n Y^m Z^k) = \left\{ \begin{array}{ll}
\frac{m!}{(m-k)!} X^{n+ m-k} & \mbox{if $m \geq k$} \\
0 & \mbox{if $k > m$ }
\end{array}
\right. 
\end{equation}
and antipode by
\begin{equation}
\tau(X^n Y^m Z^k) = (-1)^k X^m Y^n Z^k. 
\end{equation}
\end{rm}
\end{example} 

\section{Discussion}

B\"ohm and Brzezi\'nski \cite[A.1]{BB} generalize the construction $A^e \b H$
in the previous section
to a certain module algebra $A$ w.r.t.\ the action of a Hopf algebroid $H$
which is twisted by an $A$-valued cocycle on $H$.  

Panaite and Van Oystaeyen \cite{OP} observe that
the Hopf algebroid $A^e \b H$ constructed in
the last section is isomorphic to the 
Hopf algebroid  $A \odot H \odot A$
in Connes-Moscovici \cite{CM} with antipode given in
\cite{KR}, which arises in a quite different context.
The algebra $A \odot H \odot A$ formed from
a Hopf algebra $H$ and a left $H$-module algebra $A$ is linearly just $A \otimes H \otimes A$ with multiplication
given by
\begin{equation}
(a \odot h \odot b)(c \odot k \odot d) =
a (h\1 \cdot c) \odot h\2 k \odot (h\3 \cdot d)b.
\end{equation}
Note the algebra homomorphism $f_1: A \ltimes H \to
A \odot H \odot A$ given by $f_1(a \# h) := a \odot h \odot 1_A$.  Note that $f_2: H \rtimes A^{\rm op}
\to A \odot H \odot A$ given by
\begin{equation}
f_2(k \# \overline{b}) := 1_A \odot k\1 \odot k\2 \cdot b
\end{equation}
is an algebra homomorphism satisfying with $f_1$
 the hypotheses of Prop.~\ref{prop-universal}.  This leads to a
mapping $F: A^e \b H \to A \odot H \odot A$ 
given by
\begin{equation}
\label{iso-OP}
a \o \overline{b} \b h \longmapsto a \odot h\1 \odot h\2 \cdot b,
\end{equation}
which is the isomorphism in \cite[2.4]{OP}. 

Comparing the two isomorphic Hopf algebroids
(see \cite{OP} for details) 
we note that the antipode in $A^e \b H$ is given
by a simpler formula, while the $A$-$A$-bimodule structure in $A \odot H \odot A$ is
simpler. The multiplication in $A^e \b H$
is closer to the smash product of a Hopf algebra with
a bimodule algebra, which is the method of proof in \cite{OP}.

It should also be noted that \cite[3.1, 3.2]{OP} provides an equivalent condition to that in  proposition~\ref{prop-universal} which shows $A^e \b H$ is a certain universal \textit{bialgebroid}.  

Let us emphasize the picture of universals for bialgebroids over a fixed base ring
$A$.
  As observed in \cite{Lu},
for any (finite projective) algebra $A$ there
is a homomorphism of bialgebroids $A^e \to \End A$,
where $x \o \overline{y} \mapsto \lambda_x \circ \rho_y$, since $\End A$ is a terminal object
in a category of $A$-bialgebroids (existence in  \cite[Prop.\ 3.7]{Lu}, uniqueness: an easy argument).    
For similar reasons, $A^e$ is an initial object in this category.
For $A$ a left $H$-module algebra, this
homomorphism factors through the bialgebroids $A^e\b H$, $A \odot H \odot
A$, or any bialgebroid over $A$ 
as follows.    

Let $\mathcal{S}$ be a bialgebroid over $A$
with source, target mappings $s_L, t_L: A \to \mathcal{S}$
and counit $\eps: \mathcal{S} \to A$.  In addition
to Lu's mapping above, define bialgebroid
arrows $A^e \to \mathcal{S}$, $a \o \overline{b} \mapsto s_L(a) t_L(b)$ and 
the Xu anchor mapping
$\mathcal{S} \to \End A$ given by $x \mapsto
\eps(? s_L(x))$.  P.\ Xu's anchor map \cite{PX} corresponds to the
action of $\mathcal{S}$ on $A$ via source and counit \cite[3.7]{Lu},
for which $A$ becomes the unit module in the tensor category
of $\mathcal{S}$-modules. 

\begin{prop}
The natural arrows defined above form a commutative triangle of
bialgebroid homomorphisms.  
 \end{prop}
$$\begin{diagram}
A^e & &  \rTo & & \End A \\
& \SE & &  \NE  &  \\
 & & \mathcal{S}  && 
\end{diagram}$$
\begin{proof}
This follows readily from  bialgebroid identities such as 
$\eps \circ s_L = \id_A = \eps \circ t_L$.   
\end{proof}

The anchor mapping $A^e \b H \to \End A$ is given by ($a,b \in A$,
$h \in H$)
\begin{equation}
a \o \overline{b} \b h \longmapsto \lambda_a \circ \lambda_{ h \lact} \circ \rho_b
\end{equation}
where $\lambda_{ h \lact}$ denotes the endomorphism
given by left action by $h$, $x \mapsto h \lact x$. (Note that this is an algebra homomorphism since for
$a,b,c,d \in A$, $h,k \in H$, 
$$  \lambda_a \circ \lambda_{ h \lact} \circ \rho_b
\circ  \lambda_c \circ \lambda_{ k \lact} \circ \rho_d
 = \lambda_{a (h\1 \lact c)} \circ \lambda_{h\2 k\1 \lact}
\circ \rho_{d(S(k\2)\lact b)}, $$
which is the
 eq.~(\ref{eq: product}) up to a simple
re-writing.) 

 Xu's anchor mapping for  the Connes-Moscovici
bialgebroid $A \odot H \odot A$ is the mapping 
$A \odot H \odot A \to \End A$  given by
 sending $a \odot h \odot b$
into the endomorphism 
\[
\eps((a \odot h \odot b)(x \odot 1_H \odot 1_A)) = 
\eps(a(h\1 \lact x) \odot h\2 \odot b) = a(h\1 \lact x) \eps(h\2) b = \lambda_a \circ \rho_b \circ \lambda_{h \lact} (x).
\]
Note that in $\End A$ we have 
\begin{equation}
\lambda_a \circ \lambda_{h \lact} \circ \rho_b = \lambda_a \circ \rho_{h\2 \lact b} \circ \lambda_{h\1 \lact} ,
\end{equation}
which lifts to the isomorphism~(\ref{iso-OP}) $A^e \b H \rightarrow A \odot H \odot A$.

We propose a generalization of pseudo-Galois extension to \textit{pseudo-Hopf-Galois extension}
as follows. Let $H$ be a finite dimensional (or finite projective) 
Hopf algebra acting from the left on an $H$-module algebra $A$, and $B$ be a subalgebra contained in the subalgebra of invariants $$A^H = \{ b \in A \, : \,
\forall \, h \in H,\ h \lact b = \eps(h)b \}. $$
With $R := C_A(B)$ denoting the centralizer as usual, note that $H$ restricts to an action on $R$.
To be a pseudo-Hopf-Galois extension, we require
the algebra extension $A \| B$  be D2, and we require the bialgebroid homomorphism 
\begin{equation}
R^e \b H \to \End {}_BA_B, \ \ \ r \o \overline{s} \b h \mapsto \lambda_r \circ \lambda_{h \lact} \circ \rho_s
\end{equation}  
to be surjective. For example, if $A \| B$ is a Hopf-Galois extension (technically, right $H^*$-Galois), it is pseudo-Hopf-Galois since it is D2 
and by \cite{KS} 
\begin{equation}
\label{eq: nat}
\Psi:\ R \ltimes H \stackrel{\cong}{\longrightarrow} \End {}_BA_B, \ \ \ \Psi(r \# h) := \lambda_r \circ \lambda_{h \lact}.
\end{equation}
 In addition, if $A \| B$ is H-separable,
it is pseudo-Hopf-Galois since it is D2 \cite{KS}
and $R^e \cong \End {}_BA_B$ via $r \o \overline{s} \mapsto \lambda_r \circ \rho_s$.  These are the
two examples we wish to generalize at once.

The following is a third class of example of a
pseudo-Hopf-Galois extension. Let  $A \| B$ have a
split injective Galois mapping $\beta: A \o_B A \to
A \o H^*$ as $A$-$B$-bimodules and let
its trace function $A \to B$ be (a non-surjective)  Frobenius homomorphism \cite[chs.\ 4, 8]{Mo}.  Then $A \| B$ is D2 and the mapping in eq.~(\ref{eq: nat}) is a
split epimorphism via the commutative square below.

\[
\begin{diagram} A \otimes_B A & \rTo^{ \beta} & A \otimes H^* \\
 \dTo^{\cong}  & &   \dTo_{\cong}\\
  \End(A_B) & \lTo^{ \Psi} & A \# H 
\end{diagram}
     \]
\vspace{.25cm}


\end{document}